\documentclass[11pt,a4paper]{amsart}
\usepackage[english]{babel}
\usepackage[T1]{fontenc}
\usepackage[utf8]{inputenc}
\usepackage{lmodern}
\usepackage{amsmath,amsfonts,amssymb}
\usepackage[pdftex]{hyperref}
\usepackage{setspace}
\usepackage[usenames,dvipsnames]{xcolor}
\usepackage{tikz}

\title[Curvature of K\"ahler moduli]%
{Cohomological expression of\\
the curvature of K\"{a}hler moduli}
\author{Gunnar Þór Magnússon}
\thanks{Email: \texttt{gunnar@magnusson.io}}
\date{\today}

\newtheorem{theo}{Theorem}[section]
\newtheorem{prop}[theo]{Proposition}
\newtheorem{coro}[theo]{Corollary}
\newtheorem{lemm}[theo]{Lemma}
\newtheorem*{theo*}{Theorem}
\theoremstyle{definition}
\newtheorem*{defi}{Definition}

\theoremstyle{remark}
\newtheorem*{rema}{Remark}

\newcommand{\RR}{\mathbb{R}}

\newcommand{\CC}{\mathbb{C}}

\newcommand{\End}{\mathop{\mathrm{End}}}
\newcommand{\Aut}{\mathop{\mathrm{Aut}}}
\newcommand{\id}{\mathop{\mathrm{id}}}
\newcommand{\Vol}{\mathop{\mathrm{Vol}}}

\def\half{\tfrac12}
\def\onfo{\tfrac14}

\def\coho#1{\mathrm{H}^{#1}}

\def\levi{\nabla}
\def\conn{\nabla}

\def\kf{\omega}
\def\Lef{\Lambda}

\def\ton{u}
\def\ttw{v}
\def\tth{z}
\def\tfo{w}

\def\^#1{^{[#1]}}

\def\KC{C}
\def\RKC{\mathcal{\KC}}

\begin{document}

\begin{abstract}
The K\"{a}hler cone of a compact K\"ahler manifold carries a natural Riemannian
metric, given by the intersection product of its cohomology ring. We give
cohomological expressions for the Levi-Civita connection and curvature tensor of
this metric, and determine when the metric is complete.
\end{abstract}

\maketitle

\section*{Introduction}

Let $X$ be a compact K\"{a}hler manifold of complex dimension $n$. The K\"{a}hler
cone of $X$ is the set of K\"{a}hler classes, that is, $(1,1)$-classes that contain
a K\"{a}hler metric. Each K\"{a}hler class defines an inner product on the space of
$(1,1)$-classes and letting the classes vary defines a natural Riemannian metric
on the K\"{a}hler cone. 

This metric has been studied by Wilson~\cite{Wilson}, Totaro~\cite{totaro},
Wilson and Trenner \cite{WilsonTrenner} and myself~\cite{Magnusson},
sometimes by embedding the K\"{a}hler cone into the space of smooth Hermitian
metrics on the manifold via the Aubin--Calabi--Yau theorem, and sometimes by
working inside the cohomology ring.

By working with smooth forms, one can compute the curvature tensor of this
metric. Wilson~\cite{Wilson} made one such approach, using very interesting
tools that I haven't seen deployed since. I~\cite{Magnusson} made another
attempt, using the $L^2$ metric on the infinite-dimensional space of all
K\"ahler metrics. Ultimately, I cannot see that either approach descends again
to the level of cohomology (at least in a way that makes the connection to the
cohomology classes one started from obvious).

Those working in the cohomology ring have obtained more complete results.
Wilson and Trenner~\cite{WilsonTrenner} performed extensive computations in the
cohomology of Calabi--Yau threefolds. Huybrechts~\cite{Huybrechts} focused his
attention on the variation of primitive forms in the K\"{a}hler cone as the K\"{a}hler
classes vary, which is very related to the metric we're interested in and its
curvature. Totaro~\cite{totaro} also considered this metric in the more general
setting of Hessian metrics arising from homogeneous polynomials.  What unites
these efforts in the cohomology ring is their fearlessness in picking convenient
local bases to work in. Wilson, Trenner and Totaro all end up with local
coordinate expressions for the curvature tensor of the metric. Huybrechts does
not, but he doesn't compute the curvature tensor, for his sights are set
elsewhere.

The starting point of this paper was a slight feeling of dissatisfaction with
these local coordinate expressions, as I felt I didn't
understand what was going on. In
our treatment we thus stay entirely within the cohomology ring of an arbitrary
compact K\"{a}hler manifold, and avoid picking coordinates at all. This shifts
the difficulty from either previous approach from dealing with smooth forms or
complicated polynomials to computing various cup products in the cohomology
ring. Those turn out to be pleasantly manageable. 

The main novelty of this approach is that we obtain explicit and clear formulas
for the Levi-Civita connection of the metric, and for its curvature tensor. For
example, at a point $\kf$ with adjoint Lefschetz operator $\Lambda$, and for
real $(1,1)$-classes $\ton, \ttw, \tth, \tfo$, the curvature tensor is
$$
R(\ton,\ttw,\tth,\tfo)
=
- \onfo \langle \Lef(\ton \cup \tfo), \Lef(\ttw \cup \tth) \rangle
+ \onfo \langle \Lef(\ton \cup \tth), \Lef(\ttw \cup \tfo) \rangle.
$$
However, that is where the good times end. We are unable to profit from these
explicit formulas to improve on the previous general results on the curvature of
the metric in any way. The main difficulty that frustrates us is plainly visible
in the above formula: It involves the cup product in the cohomology ring of an
arbitrary compact K\"{a}hler manifold, of which next to nothing can be said. To be
able to compute the various derived tensors of the curvature tensor (the Ricci
or scalar curvatures), or to estimate their magnitudes, one needs to be able to
estimate the norm of the cup product of two classes in terms of the norms of the
individual classes uniformly over the K\"{a}hler cone. As Wilson and Trenner show
with an example, this is impossible in general; and I know of no
conditions that one can impose on the manifolds under consideration that
restrict the cohomology ring in suitable ways.

This note contains no truly new results and -- beyond the cohomological formulas
for the Levi-Civita connection and curvature tensor, which I at least find
pretty -- nothing of real interest. I have thought about this problem on and
off for around ten years, and am by now convinced there is nothing to find here.
The motivation for publishing this note is twofold: it could save other people
from wasting their time here; or if they insist, clear their path a little.

\section{The K\"{a}hler cone}

Let $X$ be a compact K\"{a}hler manifold of dimension $\dim_{\CC} X = n$.

\begin{defi}
The \emph{K\"{a}hler cone} of $X$ is the set
\begin{equation*}
\RKC(X) = \{ \kf \in \coho{1,1}(X,\RR) 
\mid
\text{$\kf$ contains a K\"{a}hler metric}
\}.
\end{equation*}
\end{defi}

If there can be no confusion about the underlying manifold $X$, we'll
just write $\RKC$ for its K\"{a}hler cone.  As the name suggests, this is
an open cones in the finite-dimensional vector space
$\coho{1,1}(X,\RR)$. The K\"{a}hler cone is the transcendental analogue
of the ample cone of a projective variety. It is described by the
transcendental version of the Nakai--Moishezon criteria due to Demailly
and Paun~\cite{DemaillyPaun}:

\begin{theo}
The K\"{a}hler cone of $X$ is a connected component of the set of real
$(1,1)$-cohomology classes that are numerically positive on analytic
cycles, that is, classes $\alpha$ such that $\int_{Z} a^p > 0$ for every
irreducible analytic set $Z$ in $X$ of dimension $p$.
\end{theo}

Since the K\"{a}hler cone $\RKC$ is an open set in a real vector space,
we can view it as a smooth manifold in its own right. It is in fact
naturally a Riemannian manifold, because each K\"{a}hler class defines an
inner product on the K\"{a}hler cone via the Hodge star operator, and these inner
products vary smoothly with the underlying class.

Let's agree on some notation before we find convenient expressions for
this metric. If $x$ is an element of the cohomology ring of $X$, we
write $x\^k := x^k/k!$ for all $k \geq 0$. This notation is
quite convenient for calculations with K\"{a}hler forms in the cohomology
ring of $X$; I learned it from Georg Schumacher.

\begin{prop}
Let $\ton,\ttw$ be elements of $T_{\kf}\RKC$. The Riemannian metric on $\RKC$
can also be defined as:
\begin{enumerate}
    \item
\hfil
$
\langle \ton, \ttw \rangle
= \Lef(\ton)\Lef(\ttw)
- \Lef\^2(\ton\wedge\ttw).
$
\hfil

    \item
The quadratic form defined by the Hessian of $-\log\Vol$.
\end{enumerate}
\end{prop}

\begin{proof}
That $(1)$ agrees with the inner product that $\kf$ defines can be seen by
taking the primitive decomposition of $\ton$ and $\ttw$, plugging it into $(1)$
and calculating until the Hodge--Riemann bilinear relations say that we have the
correct inner product.

Note that we can view $\kf$ as the tautological section $\RKC \to T_{\RKC}$
associated to the tangent bundle of any open set in a vector space.
Then $d_{\ttw} \Vol = -\Lef(\ttw)$ as we see by considering
$\Lef(\ttw) \kf\^n = \ttw \cup \kf\^{n-1}$. We then find that
$\operatorname{Hess}(\ton,\ttw) \Vol 
= d_{\ton} d_{\ttw}\Vol = \langle \ton, \ttw \rangle$ by comparing with
$(1)$.
\end{proof}

\begin{rema}
Let's write $\RKC_1 \subset \RKC$ for the set of volume-1 K\"{a}hler classes. It is a
smooth submanifold of $\RKC$, and there is a Riemannian isometry
$$
\RR \times \RKC_1 \to \RKC,
\quad
(t, \kf) \mapsto e^{t/n} \kf,
$$
where $\RR$ has the Euclidean metric and $\RKC_1$ has the restriction of the metric
on $\RKC$. As a submanifold of $\RKC$, the tangent space of $\RKC_1$ at $\kf$ is the
space of $\kf$-primitive classes.

Some authors have used this isometry to work on $\RKC_1$ when studying the metric
on the K\"{a}hler cone, as anything interesting will obviously happen there. We will
mostly leave this isometry alone and work in all of $\RKC$ instead, until the time
comes to compute the curvature tensor, when we find ourselves unable to refuse
the comforts of that subspace any longer.
\end{rema}

The theorem of Demailly and Paun describes the boundary of the
K\"{a}hler cone of a compact complex manifold. It consists of three
parts:
\begin{enumerate}
\item Limits of classes $\kf_t$ whose volume $\int_X \kf_t\^n$
tends to zero.
\item Limits of classes whose volume tends to infinity.
\item Limits of classes whose volume tends to some positive real
number, but there exists a proper irreducible complex subspace $Z
\subset X$ of dimension $p \geq 1$ whose volume tends to zero.
\end{enumerate}

Let us conspire to call $\mathcal{P} := \{\kf \in \coho{1,1}(X,\RR) \mid \kf\^n
> 0\}$ the cone of volume classes on $X$. One of its connected components
contains the K\"{a}hler cone, but is in almost all cases bigger than it.

\begin{prop}
\label{prop:fofo}
The metric on the K\"{a}hler cone of $X$ is complete if and only if the
K\"{a}hler cone is a connected component of the volume cone.
\end{prop}

\begin{proof}
We first show that the classes on the first two parts of the boundary
pose no problems. Let $I$ be an interval in the real numbers and let
$\gamma : I \to \RKC$ be a smooth path in $\RKC$ that approaches the
boundary of $\RKC$. Let $I_m = [a_m, b_m]$ be an increasing exhaustion
of $I$ by compact intervals and let $\gamma_m$ be the restriction
of $\gamma$ to $I_m$. Suppose that the volume $\Vol(X,\gamma_m)$
tends to either zero or infinity as $m$ tends to infinity.

\begin{lemm}
Let $I = [a,b]$ be a compact interval in the real numbers $\RR$,
and let $\gamma : I \to \RKC$ be a smooth path. The length of
the path $\gamma$  satisfies
$$
L(\gamma) \geq
\frac{\sqrt 2}{\sqrt n}
\bigl| \log \Vol(X,\gamma(b))
- \log \Vol(X,\gamma(a))
\bigr|.
$$
\end{lemm}

\begin{proof}[Sketch of proof.]
We apply the Cauchy--Schwarz inequality to the scalar product
$\langle\ton,\kf\rangle$; this gives
$$
|d_\ton \log \Vol(X,\kf)|^2 
= |\tfrac{1}{2}\langle\ton,\kf\rangle|^2 \leq \tfrac{n}{2} \langle\ton,\ton\rangle.
$$
Integrating and applying the triangle inequality then gives the
announced estimate.
\end{proof}

Applying the lemma on each interval $I_m$ then gives that
\begin{equation*}
L(\gamma) = \lim\limits_{m \to +\infty} L(\gamma_m) = +\infty.
\end{equation*}
Thus the limit class $\lim \gamma(t)$ on the boundary cannot be
approached by paths in $\RKC$ of finite length.

If the K\"{a}hler and volume cones of $X$ do not coincide, then there exists a
class $\alpha$ on the boundary of $\RKC$ such that $\Vol(X,\alpha) > 0$,
but there is a proper complex subspace $Z \subset X$ such that $\Vol(Z,\alpha) =
0$.

As $\alpha$ is on the boundary of the K\"{a}hler cone, then there
exists a K\"{a}hler class $\kf$ such that $\gamma(t) := \alpha +
t\kf$ is in the K\"{a}hler cone for all $t > 0$. The tangent vectors
of the path $\gamma$ are $\gamma'(t) = \kf$, and the norm of
$\gamma'(t)$ at the point $\gamma(t)$ is
$$
\displaylines{
  h(t) :=
  \langle\gamma'(t), \gamma'(t)\rangle(\gamma(t)) =
  \left(
    \frac{1}{\Vol(X,\gamma(t))}
    \int_X \kf \wedge (\alpha + t\kf)\^{n-1}
  \right)^2
\hfill\cr\hfill
{}- \frac{1}{\Vol(X,\gamma(t))}
    \int_X \kf^2 \wedge (\alpha + t\kf)\^{n-2}.
}
$$
Each of these integrals, and the function $t \mapsto \Vol(X,\gamma(t))$,
is a polynomial in $t$ on some small interval $[0,t_0]$. As $\lim_{t\to
0} \Vol(X,\gamma(t)) > 0$ the function $t \mapsto h(t)$ is continuous
and positive on a compact interval, so the integral $L(\gamma)$ of its
square root exists and is finite.
\end{proof}

A holomorphic map $f : X \to Y$ between compact K\"{a}hler manifolds induces
a morphism $f^* : \coho{*}(Y,\RR) \to \coho{*}(X,\RR)$ in cohomology
that respects the Hodge decomposition. However, if $\kf$ is a K\"{a}hler
class on $Y$, then $f^*\kf$ is hardly ever a K\"{a}hler class on $X$. This
happens mostly if $f$ is either an embedding or a finite covering map.

\begin{prop}
Let $f : X \to Y$ be a finite surjective morphism. Let $g_X$ and $g_Y$
be the metrics on the K\"{a}hler cones of $X$ and $Y$, respectively.
Then the pullback morphism $f^* : \RKC(Y) \to \RKC(X)$ is a Riemannian
embedding.
\end{prop}

\begin{proof}
Let $\kf$ be a point in $\RKC(Y)$. The volume of $X$ with respect to $f^*\kf$ is
\begin{equation*}
  \Vol(X,f^*\kf) = p \, \Vol(Y,\kf)
\end{equation*}
as $f$ is finite of degree $p$. It follows that $f^*$ is an embedding.
\end{proof}

\begin{coro}
The group $\Aut X$ of holomorphic automorphisms of $X$ acts by
isometries on the K\"{a}hler cone $\RKC(X)$.
\end{coro}

A closer look reveals that this last statement contains less information
than first meets the eye. The automorphism group $\Aut X$ of a compact
complex manifold is a Lie group and it splits roughly into two parts; a
positive-dimensional group given by the flows of holomorphic vector
fields, or elements of $\coho{0}(X,T_X)$, and a discrete part consisting of
``other'' automorphisms. The isomorphisms generated by vector fields act
trivially on the cohomology ring of $X$, so the only part of $\Aut X$
that possibly acts by nontrivial isometries on $\RKC(X)$ is discrete.

\section{Connection and curvature}

We start with a couple of preliminary computations.

\begin{lemm}
If $u_1, \ldots, u_k$ are real $(1,1)$-classes, then
$$
\displaylines{
d_{\ttw} \Lef\^{k}(u_1 \cup \cdots \cup u_k)
= -\Lef(\ttw) \Lef\^{k}(u_1 \cup \cdots \cup u_k)
\hfill\cr\hfill{}
+ \Lef\^{k}(d_{\ttw} u_1 \cup \cdots \cup u_k)
+ \cdots +
\Lef\^{k}(u_1 \cup \cdots \cup d_{\ttw} u_k)
\hfill\cr\hfill{}
+ \Lef\^{k+1}(u_1 \cup \cdots \cup u_k \cup \ttw).
}
$$
\end{lemm}

\begin{proof}
This is clear once we write
$$
\Lef\^{k}(u_1 \cup \cdots \cup u_k)
= \frac{1}{\Vol(X,\kf)} \int_X u_1 \cup \cdots \cup u_k \cup \kf\^{n-k}
$$
and compute.
\end{proof}

\begin{lemm}
\label{lemm:triple}
Let $\ton,\ttw,\tth$ be $(1,1)$-classes. Then
$$
\langle \Lef(\ton\cup\ttw), \tth \rangle
= - \Lef\^3(\ton \cup \ttw \cup \tth)
+ \Lef\^2(\ton \cup \ttw) \Lef(\tth).
$$
\end{lemm}

\begin{proof}
First note that if $\tth$ is a $(1,1)$-class, then $\tth = (\tth -
\frac1n \Lef(\tth)\kf) + \frac1n \Lef(\tth) \kf$ is its primitive
decomposition. Then
\begin{align*}
*(\kf \cup \tth)
&= *( \kf \cup (\tth - \tfrac1n \Lef(\tth)\kf))
+ * \Bigl( \tfrac{2}{n} \Lef(\tth) \kf\^2 \bigr)
\\
&= -(\tth - \tfrac1n \Lef(\tth)\kf) \cup \kf\^{n-3}
+ \tfrac{2}{n} \Lef(\tth) \kf\^{n-2}
\\
&= -\tth \cup \kf\^{n-3}
+ \tfrac{n-2}{n} \Lef(\tth) \kf\^{n-2}
+ \tfrac{2}{n} \Lef(\tth) \kf\^{n-2}
\\
&= -\tth \cup \kf\^{n-3}
+ \Lef(\tth)\, \kf\^{n-2}.
\end{align*}
We now get
\begin{align*}
\langle \Lef(\ton \cup \ttw), \tth \rangle \,\kf\^n
&= \Lef(\ton \cup \ttw) \cup 
(-\tth \cup \kf\^{n-2} + \Lef(\tth)\, \kf\^{n-1})
\\
&= - \Lef\^2(\Lef(\ton \cup \ttw) \cup \tth) \,\kf\^{n}
+ 2 \Lef\^2(\ton \cup \ttw) \Lef(\tth) \,\kf\^{n},
\end{align*}
which proves the result.
\end{proof}

Recall that the Levi-Civita connection is the unique connection on the tangent
bundle that's compatible with the metric and is torsion-free. That is, it
satisfies
$$
d \langle\ton, \ttw\rangle 
= \langle\levi \ton, \ttw\rangle + \langle\ton, \levi \ttw\rangle,
\quad
\levi_{\ton}\ttw - \levi_{\ttw}\ton = [\ton,\ttw]
$$
for all sections $\ton, \ttw$ of the bundle.

\begin{prop}
\label{prop:connection}
The Levi-Civita connection of the Riemannian metric $g$ on $\RKC$ is
$$
\levi_{\tth} \ton
=
d_\tth \ton
-\tfrac12 \Lef(\ton) \tth
-\tfrac12 \Lef(\tth) \ton
+\tfrac12 \Lef(\ton \cup \tth).
$$
\end{prop}

\begin{proof}
The connection we've written down satisfies $\conn_{\ton} \ttw - \conn_{\ttw}
\ton = [\ton, \ttw]$ by inspection. We turn to its computation.

The metric is defined by
$$
\langle \ton, \ttw \rangle
= \Lef \ton \, \Lef \ttw
- \Lef\^2(\ton \cup \ttw).
$$
Taking the derivative of this in the $\tth$ direction gives
$$
\displaylines{
d_{\tth} \langle \ton, \ttw \rangle
= 
- \Lef(\ton) \Lef(\ttw) \Lef(\tth)
+ \Lef(d_{\tth}\ton) \Lef(\ttw)
+ \Lef\^{2}(\ton \cup \tth) \Lef(\ttw)
\hfill\cr\hfill{}
- \Lef(\ton) \Lef(\ttw) \Lef(\tth)
+ \Lef(\ton) \Lef(d_{\tth}\ttw)
+ \Lef(\ton) \Lef\^{2}(\ttw \cup \tth)
\cr\hfill{}
+ \Lef(\tth) \Lef\^2(\ton \cup \ttw)
- \Lef\^2(d_{\tth} \ton \cup \ttw)
- \Lef\^2(\ton \cup d_{\tth} \ttw)
- \Lef\^3(\ton \cup \ttw \cup \tth)
\cr{}
\phantom{d_{\tth} \langle \ton, \ttw \rangle}
=: \langle d_{\tth} \ton, \ttw \rangle
+ \langle \ton, d_{\tth} \ttw \rangle
+ A(\ton, \ttw, \tth).
\hfill
}
$$
We're going to write $A = \frac12 A + \frac12 A$ and try to write the first half
as an inner product with $\ttw$, and the second half as an inner product with
$\ton$.

To that end, we note that
$$
\displaylines{
A(\ton, \ttw, \tth)
= 
- \Lef(\ton) \Lef(\ttw) \Lef(\tth)
+ \Lef\^{2}(\ton \cup \tth) \Lef(\ttw)
\hfill\cr\hfill{}
- \Lef(\ton) \Lef(\ttw) \Lef(\tth)
+ \Lef(\ton) \Lef\^{2}(\ttw \cup \tth)
\cr\hfill{}
+ \Lef(\tth) \Lef\^2(\ton \cup \ttw)
- \Lef\^3(\ton \cup \ttw \cup \tth)
\cr{}
\phantom{A(\ton, \ttw, \tth)}
= 
- \Lef(\tth) \langle \ton, \ttw \rangle
- \Lef(\ton) \langle \tth, \ttw \rangle
+ \Lef\^{2}(\ton \cup \tth) \Lef(\ttw)
- \Lef\^3(\ton \cup \ttw \cup \tth)
\hfill\cr\hfill{}
=
- \langle \Lef(\tth) \ton, \ttw \rangle
- \langle \Lef(\ton) \tth, \ttw \rangle
+ \langle \Lef(\ton \cup \tth), \ttw \rangle
}
$$
by Lemma~\ref{lemm:triple}, which can indeed be written as an inner product with
$\ttw$. Since $A$ is symmetric in $\ton, \ttw, \tth$, we can start again from
$A$ and write it as an inner product with $\ton$. Taking half of each, we arrive
at our claimed form of the connection.
\end{proof}

\begin{coro}
\label{coro:kahlerform}
\begin{itemize}
    \item 
$\levi \kf = 0$.
    \item 
If $\ton$ is a primitive vector field, then $\levi \ton$ is also
primitive.
\end{itemize}
\end{coro}

\begin{theo}
\label{theo:curvature}
The curvature tensor of the metric on the K\"{a}hler cone is
\begin{equation*}
R(\ton,\ttw,\tth,\tfo)
= 
- \onfo \langle \Lef(\ton \cup \tfo), \Lef(\ttw \cup \tth) \rangle
+ \onfo \langle \Lef(\ton \cup \tth), \Lef(\ttw \cup \tfo) \rangle.
\end{equation*}
\end{theo}

\begin{proof}
We may assume that all the tangent fields $\ton,\ttw,\tth,\tfo$ are
primitive, either by appealing to the isometric splitting of the K\"{a}hler
cone, by using that $\conn \kf = 0$ and the symmetries of the curvature
tensor to see that $R$ always degenerates to the primitive parts of our
classes, or by calculating the curvature tensor first for primitive
classes and then doing painful algebra to see that the general case
degenerates to that one. However we do it, we find that for primitive
fields we have
$$
\conn_{\ttw} \tth
= d_{\ttw} \tth + \half \Lef(\ttw \cup \tth)
$$
and
$$
\displaylines{
\conn_{\ton} \conn_{\ttw} \tth
= d_{\ton} d_{\ttw} \tth 
+ \half (d_{\ton} \Lef)(\ttw \cup \tth)
+ \half \Lef(d_{\ton} \ttw \cup \tth)
\hfill\cr\hfill{}
+ \half \Lef(\ttw \cup d_{\ton} \tth)
+ \half \Lef(\ton \cup d_{\ttw} \tth)
+ \onfo \Lef(\ton \cup \Lef(\ttw \cup \tth)).
}
$$
This gives
$$
R(\ton,\ttw) \tth
= 
\half (d_{\ton} \Lef)(\ttw \cup \tth)
+ \onfo \Lef(\ton \cup \Lef(\ttw \cup \tth))
- \half (d_{\ttw} \Lef)(\ton \cup \tth)
- \onfo \Lef(\ttw \cup \Lef(\ton \cup \tth))
$$
since the other terms either make up $\conn_{[\ton,\ttw]}\tth$ or are
symmetric in $\ton,\ttw$. To make sense of this, it's convenient to take
the inner product with $\tfo$.

We have $*(\kf \cup \tfo) = -\kf\^{n-3} \cup \tfo$ since $\tfo$ is
primitive, so 
$$
\langle \Lef(\ttw \cup \tth), \tfo \rangle
= \langle \ttw \cup \tth, \kf \cup \tfo \rangle
= -\Lef\^3(\ttw \cup \tth \cup \tfo).
$$
Differentiating this in the direction of $\ton$ gives
$$
\displaylines{
\langle (d_{\ton}\Lef)(\ttw \cup \tth), \tfo \rangle
+ \half \langle \Lef(\ton \cup \Lef(\ttw \cup \tth)), \tfo \rangle
+ \half \langle \Lef(\ttw \cup \tth), \Lef(\ton \cup \tfo) \rangle
\hfill\cr\hfill{}
= 
-\Lef\^4(\ton \cup \ttw \cup \tth \cup \tfo)
}
$$
after canceling out the terms that involve the derivatives of the
tangent fields. Then one part of the curvature tensor is
$$
\displaylines{
\half \langle (d_{\ton}\Lef)(\ttw \cup \tth), \tfo \rangle
+ \onfo \langle \Lef(\ton \cup \Lef(\ttw \cup \tth)), \tfo \rangle
\hfill\cr\hfill{}
=- \half \Lef\^4(\ton \cup \ttw \cup \tth \cup \tfo)
- \onfo \langle \Lef(\ttw \cup \tth), \Lef(\ton \cup \tfo) \rangle.
}
$$
The first term is symmetric in $\ton,\ttw$, so we get
$$
R(\ton,\ttw,\tth,\tfo)
= 
- \onfo \langle \Lef(\ton \cup \tfo), \Lef(\ttw \cup \tth) \rangle
+ \onfo \langle \Lef(\ton \cup \tth), \Lef(\ttw \cup \tfo) \rangle
$$
as promised.
\end{proof}

\begin{rema}
If $x,y$ are $(2,2)$-classes, then 
$$
\Lef\^4(x \cup y) 
= \langle x, y \rangle 
- \langle \Lef(x), \Lef(y) \rangle 
+ \langle \Lef\^2(x), \Lef\^2(y) \rangle;
$$
see \cite{magnusson_inner_product}.
An alternate expression for the curvature tensor is thus
\begin{equation*}
\displaylines{
R(\ton,\ttw,\tth,\tfo)
= 
- \onfo \langle \ton, \tfo \rangle 
    \langle \ttw, \tth \rangle
+ \onfo \langle \ton, \tth \rangle 
    \langle \ttw, \tfo \rangle
\hfill\cr\hfill{}
- \onfo \langle \ton \cup \tfo, \ttw \cup \tth \rangle
+ \onfo \langle \ton \cup \tth, \ttw \cup \tfo \rangle,
}
\end{equation*}
so the curvature tensor is a perturbation of the curvature tensor of a space
form of constant sectional curvature. Unfortunately there is no known way
to control the perturbation terms in general; at least bounding them from above is
impossible by Wilson and Trenner's example~\cite{WilsonTrenner}.
\end{rema}

\subsection*{On algebraic curvature tensors}

The expression for the curvature tensor suggests that we could investigate the
operation $(\ton, \ttw) \mapsto \half \Lef(\ton \cup \ttw)$ to understand the
curvature of the metric. This operation defines an algebra structure on
$\coho{1,1}(X,\RR)$; as an algebra, it is commutative, non-associative and
non-unital (if it had a unit, it would have to be a multiple of $\kf$, which
doesn't work). This algebra structure varies as the K\"{a}hler class $\kf$ varies.

This curvature tensor conforms to a form of algebraic curvature tensors that, as
far as I know, have not received much attention. We gather here some
trivialities about them, the first of which suggests this will not be a fertile
line of investigation.

\begin{prop}
Let $V$ be a real vector space, equipped with an inner product $\langle \ ,\,\rangle$, and an algebra structure $(x,y) \mapsto x\cdot y$. If the algebra structure is commutative, then
$$
R(x,y,z,w) := \langle x \cdot w, y \cdot z \rangle - \langle x \cdot z, y \cdot w \rangle
$$
is an algebraic curvature tensor.
\end{prop}

\begin{proof}
  It is immediate that $R(y,x,z,w) = R(x,y,w,z) = -R(x,y,z,w)$. The commutativity entails that $R(z,w,x,y) = R(x,y,z,w)$, so $R$ defines a symmetric bilinear form on $\bigwedge^2V$. The commutativity also entails that $R$ satisfies the Bianchi identity
  $$
  \displaylines{
    \phantom{\qedsymbol}
    \hfill
  R(x,y,z,w) + R(y,z,x,w) + R(z,x,y,w) = 0.
  \hfill
  \qedhere
  }
  $$
\end{proof}

The moral of this proposition is perhaps that one should not expect to be able to prove very much about our curvature tensor from formal properties alone. After all, the commutative algebra structures on a vector space of dimension $h^{1,1}$ form a vector space of dimension $(h^{1,1})^2(h^{1,1}+1)/2$, while we can at best expect to generate a space of dimension $h^{1,1}$ therein by deforming our K\"ahler classes. Without any way of distinguishing the structures defined by a K\"ahler class (if any) from the others, and without this class of curvature tensors having some special properties, there is then little hope of progress in this direction.

Recall that if $b$ is a symmetric bilinear form on $V$, then the Kulkarni--Nomizu product of $b$ is the algebraic curvature tensor defined as
$$
(b \wedge b)(x,y,z,w)
= b(x,z)b(y,w) - b(x,w)b(y,z).
$$

\begin{prop}
A curvature tensor defined by an algebra is a sum of Kulkarni--Nomizu products.
\end{prop}

\begin{proof}
Let $(x_1,\ldots,x_n)$ be an orthonormal basis of $V$. Define bilinear forms $b_l(x,y) := \langle x \cdot y, x_l \rangle$. These are symmetric as the algebra is symmetric, and satisfy $x \cdot y = \sum_{l=1}^n b_l(x,y) x_l$. It follows that
$$
\displaylines{
  \phantom{\qedsymbol}
  \hfill
  R(x,y,z,w) = -\sum_{l=1}^n (b_l \wedge b_l)(x,y,z,w).
  \hfill
  \qedhere
}
$$
\end{proof}

The curvature tensors of our algebra structures are made up of symmetric bilinear forms on $S^2V$.
If $x, y$ are vectors in $V$, we'll write $xy$ for the induced vector $\frac12(x \otimes y + y \otimes x)$ in $S^2V$. The inner product on $V$ induces an inner product on $S^2V$ by
$$
\langle xy, zw \rangle
=
\tfrac12
\bigl(
\langle x, z \rangle \langle y, w \rangle
+ \langle x, w \rangle \langle y, z \rangle
\bigr).
$$
We'll say that an algebraic curvature tensor has constant sectional curvature if it is equal to a multiple of the Kulkarni--Nomizu product of an inner product.

\begin{prop}
An algebraic curvature tensor induced by an algebra structure has constant sectional curvature if and only if there exists a scalar $\lambda$ and a symmetric $4$-tensor $A$ such that
$$
\langle x \cdot y, z \cdot w \rangle
= 2\lambda \langle xy, zw \rangle + A(x,y,z,w)
$$
for all vectors $x,y,z,w \in V$.
\end{prop}

\begin{proof}
  If this condition holds, it is a simple computation to show that the curvature tensor defined by the algebra has constant sectional curvature $-\lambda$ with our sign choices.

  Conversely, if the curvature tensor has constant sectional curvature $-\lambda$, the symmetries of the curvature tensor entail that the linear form
$$
A(x,y,z,w) := \langle x \cdot w, y \cdot z \rangle - 2\lambda \langle xw, yz \rangle
$$
is symmetric.
\end{proof}

Recall that a derivation $D$ of an algebra is a linear map on the underlying vector space such that $D(x \cdot y) = Dx \cdot y + x \cdot Dy$ for all vectors $x,y$. The derivations of an algebra on $V$ form a subalgebra of $\End V$. One can constrain them a little by formal manipulations in the case of our algebra:

\begin{prop}
If $D$ is a derivation of the algebra on $\coho{1,1}(X,\mathbb{R})$, then $D\omega = 0$, $Dx$ is primitive for all $x$, and $D^t = -D$.
\end{prop}

\begin{proof}
We have $[L, \Lambda] = (2-n)\id$ on $(1,1)$-forms. We can interpret this as a statement about the algebra product of a class with $\omega$:
$$
x \cdot \omega
= \tfrac12 \Lambda(x) \omega + \tfrac12 (n-2) x.
$$
We have $\omega \cdot \omega = (n-1) \omega$. Then
$$
(n-1) D\omega
= 2 D\omega \cdot \omega
= \Lambda(D\omega) \, \omega + n D\omega.
$$
Then
$$
D\omega = -\Lambda(D\omega) \, \omega,
$$
that is, $D\omega$ is a multiple of $\omega$.

Suppose that $u$ is primitive. Then $u \cdot \omega = \tfrac12 (n-2) u$. We get
\begin{align*}
\tfrac12 (n-2) Du
&= Du \cdot \omega + u \cdot D\omega
\\
&= \tfrac12 \Lambda(Du) \, \omega 
+ \tfrac12 (n-2) Du
 + \Lambda(D\omega) \tfrac12 (n-2) u,
\end{align*}
That is,
$$
\Lambda(Du) \, \omega + (n-2) \Lambda(D\omega) \, u = 0.
$$
As $\omega$ and $u$ are orthogonal, the only way this can hold is if $\Lambda(Du) = \Lambda(D\omega) = 0$. Then we also get $D\omega = 0$.

It follows that $Dx$ is primitive for any $(1,1)$-class $x$. For primitive classes $u$ and $v$, we have
$$
D(u \cdot v) = Du \cdot v + u \cdot Dv.
$$
Note that if either $u$ or $v$ is a primitive class, then $\Lambda(u \cdot v) = -\langle u, v \rangle$. From the above, it follows that
$$
\langle Du, v \rangle + \langle u, Dv \rangle = 0,
$$
so the linear morphism $D$ satisfies $D^t = -D$.
\end{proof}

\bibliographystyle{alpha}
\bibliography{main}

\end{document}